\documentclass[10.5pt]{amsart}
\usepackage{amsmath}
\usepackage[applemac]{inputenc}
\usepackage[frenchb]{babel}
\usepackage{amscd}
\usepackage{xspace}
\usepackage{verbatim}
\newcommand{\fnz}{\footnotesize}

\newcommand{\be}{\begin{equation}}
\newcommand{\bel}[1]{\begin{equation}\label{#1}}
\newcommand{\ee}{\end{equation}}%% This macro does not work with amstex.
\newtheorem{subn}{\name}

\newcommand{\bsn}[1]{\def\name{#1}\begin{subn}}
\newcommand{\esn}{\end{subn}}
\newtheorem{sub}{\name}%[section]
\newcommand{\bs}{\begin{sub}}
\newcommand{\es}{\end{sub}}
\newcommand{\bsl}[1]{\begin{sub}\label{#1}}
\newcommand{\bth}[1]{\def\name{Theorem}\begin{sub}\label{t:#1}}
\newcommand{\blemma}[1]{\def\name{Lemma}\begin{sub}\label{l:#1}}
\newcommand{\bcor}[1]{\def\name{Corollary}\begin{sub}\label{c:#1}}
\newcommand{\bdef}[1]{\def\name{Definition}\begin{sub}\label{d:#1}}
\newcommand{\bprop}[1]{\def\name{Proposition}\begin{sub}\label{p:#1}}
\newcommand{\rlemma}[1]{Lemma~\ref{l:#1}}

\newcommand{\BA}{\begin{array}}
\newcommand{\EA}{\end{array}}
\newcommand{\BAN}{\renewcommand{\arraystretch}{1.2}
\setlength{\arraycolsep}{2pt}\begin{array}}
\newcommand{\BAV}[2]{\renewcommand{\arraystretch}{#1}
\setlength{\arraycolsep}{#2}\begin{array}}
\newcommand{\BSA}{\begin{subarray}}
\newcommand{\ESA}{\end{subarray}}
\newcommand{\BAL}{\begin{aligned}}
\newcommand{\EAL}{\end{aligned}}
\newcommand{\BALG}{\begin{alignat}}
\newcommand{\EALG}{\end{alignat}}%% the abbrev. does not work with latex2e
\newcommand{\BALGN}{\begin{alignat*}}
\newcommand{\EALGN}{\end{alignat*}}%% the abbrev. does not work with latex2e
\newcommand{\abs}[1]{\left |#1\right |}%% adjustable vertical delimiters
\newcommand{\opname}[1]{\mbox{\rm #1}\,}

\newcommand{\dist}{\opname{dist}}
\newcommand{\myfrac}[2]{{\displaystyle \frac{#1}{#2} }}

\newcommand{\prt}{\partial}

\newcommand{\ti}{\times}

\def\gth{\theta}                         
           
            \def\gl{\lambda}
\def\gm{\mu}

     \def\Gd{\Delta}      
\def\Gth{\Theta}
          
\def\Gw{\Omega}              

\def\CA{{\mathcal A}}

      \def\CL{{\mathcal L}}

   \def\BBR {\mathbb R}

\begin{document}%

%%Setting up the TITLE and AUTHOR
\noindent Partial Differential Equations/Differential Geometry\smallskip

\title[Quasilinear equations on manifolds]{Quasilinear elliptic Hamilton-Jacobi equations on complete  manifolds}
\maketitle
\noindent {Marie-Fran\c{c}oise Bidaut-V\'eron
\footnote{\noindent Laboratoire de Math\'{e}matiques et Physique Th\'{e}orique, CNRS UMR 7350,
Facult\'{e} des Sciences, 37200 Tours France. E-mail: veronmf@univ-tours.fr},
 Marta Garcia-Huidobro
 \footnote{\noindent 
Departamento de Matematicas, Pontifica Universidad Catolica de Chile
Casilla 307, Correo 2, Santiago de Chile. E-mail: mgarcia@mat.puc.cl}, 
Laurent V\'{e}ron
 \footnote{\noindent 
Laboratoire de Math\'{e}matiques et Physique Th\'{e}orique, CNRS UMR 7350,
Facult\'{e} des Sciences, 37200 Tours France. E-mail: veronl@univ-tours.fr}}

\date{}
%%%%%%%%%%
\begin{abstract} Let $(M^n,g)$ be a $n$-dimensional complete, non-compact and connected Riemannian manifold, with Ricci tensor $Ricc_g$ and sectional curvature $Sec_g$. Assume $Ricc_g\geq (1-n)B^2$, and either $p>2$ and $Sec_g(x)=o(dist^2(x,a))$ when $dist^2(x,a)\to\infty$ for $a\in M$, or $1<p<2$ and $Sec_g(x)\leq 0$. If $q>p-1> 0$, any $C^1$ solution of (E) $-\Gd_pu+\abs{\nabla u}^q=0$ on $M$ satisfies $\abs{\nabla u(x)}\leq c_{n,p,q}B^{\frac{1}{q+1-p}}$ for some constant $c_{n,p,q}>0$. As a consequence there exists $c_{n,p}>0$ such that any positive $p$-harmonic function $v$ on $M$ satisfies $v(a)e^{-c_{n,p}B\dist (x,a)}\leq v(x)\leq v(a)e^{c_{n,p}B\dist (x,a)}$ for any $(a,x)\in M\ti M$.
  \end{abstract}
%%%%%%%%%%%%%%%%
%\maketitle
%%%%%%%%%%%%%
\begin{center}{\bf  Equations de Hamilton-Jacobi quasilin\'{e}aires  sur une vari\'{e}t\'{e} compl\`{e}te}\end{center}
\smallskip
\begin{quotation}
 {\fnz{\sc R\'esum\'e.} Soit $(M^n,g)$ une variét\'{e} riemannienne $n$-dimensionnelle compl\`ete, non compacte et connexe de courbures de Ricci $Ricc_g$ et  sectionnelle $Sec_g$. On suppose $Ricc_g\geq (1-n)B^2$ et $Sec_g(x)=o(dist^2(x,a))$ si $dist^2(x,a)\to\infty$ pour $a\in M$ si $p> 2$, ou $Sec_g(x)\leq 0$ si $1<p<2$. Si $q>p-1>0$, toute solution de classe $C^1$ de (E) $-\Gd_pu+\abs{\nabla u}^q=0$ sur $M$ satisfait \`a $\abs{\nabla u(x)}\leq c_{n,p,q}B^{\frac{1}{q+1-p}}$ où $c_{n,p,q}>0$ est une constante. On en d\'{e}duit qu'il existe $c_{n,p}>0$ tel que toute fonction $p$-harmonique positive  $v$ sur $M$ satisfait à l'encadrement suivant, $v(a)e^{-c_{n,p}B\dist (x,a)}\leq v(x)\leq v(a)e^{c_{n,p}B\dist (x,a)}$ pour tout $(a,x)\in M\ti M$.
}
\end{quotation}
 \bigskip
%\numberwithin{equation}{section}
 {\bf Version fran\c{c}aise abr\'eg\'ee.}
 %% %%
 Soit $(M^n,g)$ une variété riemannienne complète, non-compacte et connexe de courbure de Ricci $Ricc_g$ et courbure sectionnelle 
 $Sec_g$. Pour tout $p>1$, on dénote par $u\mapsto \Gd_pu:=div\left(\abs{\nabla u}^{p-2}\nabla u\right)$ le $p$-Laplacien sur $M$ pour la métrique $g$. Notre résultat principal est le suivant\medskip
 
 \noindent{\bf Theorème 1. }{\it Soit $B\geq 0$ tel que $Ricc_g\geq (1-n)B^2$ et $q>p-1> 0$. On suppose 
  \bel{equ0}
 \lim_{dist(x,a)\to\infty}\frac{Sec_g(x)}{(dist(x,a))^2}=0
 \ee
pour tout $a\in M$ si $p>2$, ou $Sec_g\leq 0$ si $1<p<2$. Il existe alors $c_{n,p,q}>0$ telle que toute solution $u\in C^{1}(M)$ de 
 \bel{equ1}
 -\Gd_pu+\abs{\nabla u}^q=0\qquad\text{sur }M
 \ee
 vérifie
  \bel{equ2}
 \abs{\nabla u(x)}\leq c_{n,p,q}B^{\frac{1}{q+1-p}}\qquad\forall x\in M.
 \ee
}
%\begin {center}
Une des conséquences est un théorème de type Liouville. \medskip

 \noindent{\bf Corollaire 2. }{\it Supposons que $Ricc_g\geq 0$, $q>p-1> 0$ et que les hypothèses  du Th\'eor\`eme 1 portant sur la courbure sectionnelle soient vérifiées si $p\neq 2$. Alors toute solution $u\in C^{1}(M)$ de (\ref{equ1}) est constante.}\medskip
 
Si $v$ est une fonction $p$-harmonique positive sur $M$, la fonction $u:=-(p-1)\ln v$ vérifie
  \bel{equ3}
 -\Gd_pu+\abs{\nabla u}^p=0\qquad\text{sur }M.
 \ee
 En utilisant le résultat du théorème 1, on en déduit\medskip
 
 \noindent{\bf Théorème 3. }{\it Supposons que $p>1$ et que les hypothèses du Th\'eor\`eme 1 portant sur la courbure soient vérifiées. Il existe alors une constante $c_{n,p}>0$ telle que toute fonction $p$-harmonique et positive $v$ sur $M$ vérifie
 \bel{equ4}
 v(a)e^{-c(n,p)B\dist (x,a)}\leq v(x)\leq v(a)e^{c(n,p)B\dist (x,a)}\qquad\forall (a,x)\in M\ti M.
  \ee}

Quand $p=2$ Cheng et Yau \cite{CY} ont montré que toute fonction harmonique positive sur une variété riemannienne complète à courbure de Ricci positive est une constante. Dans le cas des fonctions $p$-harmoniques positives et sous l'hypothèse de minoration uniforme de la courbure sectionnelle, $Sec_g\geq -B^2$, Kotschwar et Ni \cite{KN}  montrent que  toute fonction $p$-harmonique positive $v$ sur $M$  vérifie l'estimation suivante,
 \bel{equ4'}\frac{\abs{\nabla v}}{v} \leq  (p-1)B.\ee
 Notons que leur hypothèse implique $Ricc_g\geq (1-n)B^2$.
\begin {center}
---------------------------------------------------------------------------------------------
\end {center} \setcounter{equation}{0}

Let  $(M^n,g)$ be a complete, connected and non compact Riemannian manifold with Ricci curvature $Ricc_g$ and sectionnal  
curvature $Sec_g$. For $p>1$ we denote by  $\Gd_p$ the p-Laplacian defined in the metric $g$ by
$$\Gd_pu:=div\left(\abs{\nabla u}^{p-2}\nabla u\right),$$
and thus $\Gd_2$ is the Laplace-Beltrami operator on M. If $p=2$ a classical result due to Cheng and Yau \cite{CY} asserts that if $Ricc_g$ is nonnegative, any nonnegative harmonic function $v$ is a constant. In \cite{KN}, Kotschwar et Ni obtained sharper results dealing with positive $p$-harmonic functions under the assumption that $Sec_g\geq -B^2$. They proved that if $v$ is such a function, it satisfies
 \bel{equ4-2}\frac{\abs{\nabla v}}{v} \leq  (p-1)B.\ee
Their assumption on $Sec_g$ implies $Ricc_g\geq (1-n)B^2$. They also noticed that if $p=2$ their estimate holds under  the previous lower estimate on the Ricci curvature.  In this note we give an extension of their result in imbedding it  the more general class of quasilinear Hamilton-Jacobi type equations
 \bel{equ1-2}
 -\Gd_pu+\abs{\nabla u}^q=0\qquad\text{on }M.
 \ee
Our main result is the following\medskip

 \noindent{\bf Theorem 1. }{\it Let $B\geq 0$ such that $Ricc_g\geq (1-n)B^2$. If $p>2$ we assume that for any $a\in M$
  \bel{equ0}
 \lim_{dist(x,a)\to\infty}\frac{Sec_g(x)}{(dist(x,a))^2}=0,
 \ee
and if $1<p<2$ that $Sec_g\leq 0$.  Then there exists $c_{n,p,q}>0$ such that any solution $u\in C^{1}(M)$ of  (\ref{equ1-2})
satisfies
  \bel{equ2-2}
 \abs{\nabla u(x)}\leq c_{n,p,q}B^{\frac{1}{q+1-p}}\qquad\forall x\in M.
 \ee
}\medskip

A clear consequence of (\ref{equ2}) is the following Liouville theorem\medskip

\noindent{\bf Corollary 2. }{\it Assume $Ricc_g\geq 0$ and that the assumptions of  Theorem 1 concerning  $Sec_g$  hold if $p\neq  2$. Then any solution $u\in C^{1}(M)$ of (\ref{equ1-2}) is constant.}\medskip
 
If $v$ is a positive  $p$-harmonic function on $M$, then $u:=-(p-1)\ln v$ satisfies
  \bel{equ3-2}
 -\Gd_pu+\abs{\nabla u}^p=0\qquad\text{on }M.
 \ee
Therefore estimate (\ref{equ2-2}) yields to the following result\medskip

 \noindent{\bf Theorem 3. }{\it Assume that $p>1$ and the curvature assumptions of Theorem 1 are fulfilled. Then there exists a constant $c_{n,p}>0$ such that any positive $p$-harmonic function $v$ on $M$ satisfies
 \bel{equ4'-2}
 v(a)e^{-c(n,p)B\dist (x,a)}\leq v(x)\leq v(a)e^{c(n,p)B\dist (x,a)}\qquad \forall (a,x)\in M\ti M.
 \ee}

\noindent{\bf Proof of Theorem 1.} Let $M_+:=\{x\in M:\abs{\nabla u(x)}>0\}$. Then $M_+$ is open and $u\in C^3(M_+)$ since the equation is no longer degenerate.  The proof is based upon the fact that $z=\abs{\nabla u}^2$ is a subsolution of an elliptic differential inequality with a superlinear absorption term (see \cite{PV1} for other applications). We denote by $TM$ the tangent bundle of $M$ and by $\langle.,.\rangle$ the scalar product induced by the metric $g$. We recall that any  $C^3$-function $u$ verifies the B\"ochner-Weitzenb\"ock formula; combined with Schwarz inequality it yields to
 \bel{W1}\BA {l}
\myfrac{1}{2}\Gd_2\abs{\nabla u}^2=\abs{D^2u}^2+\langle\nabla \Gd_2 u,\nabla u\rangle+Ricc_g(\nabla u,\nabla u)
\\\phantom{\myfrac{1}{2}\Gd_2\abs{\nabla u}^2}
\geq \myfrac{1}{n}\abs{\Gd_2u}^2+\langle\nabla \Gd_2 u,\nabla u\rangle+Ricc_g(\nabla u,\nabla u),
\EA\ee
where $D^2u$ is the Hessian. If $u$ is a $C^1$ solution of (\ref{equ1-2}), then  $z=\abs{\nabla u}^2$ satisfies
 \bel{W2}-\Gd_2u-\myfrac{p-2}{2}\myfrac{\langle\nabla z,\nabla u\rangle}{z}+z^{\frac{q+2-p}{2}}=0 \ee
on $M_+$. Replacing $\Gd_2u$ in (\ref{W1}) it follows that, for any $a>0$, 
\bel{W3}\BA {l}
\Gd_2 z+(p-2)\myfrac{\langle D^2z(\nabla u),\nabla u\rangle}{z}\geq\myfrac{2a^2}{N}z^{q+2-p}-\myfrac{1}{Na^2}\myfrac{\langle\nabla z,\nabla u\rangle^2}{z^2}-\myfrac{(p-2)}{2}\myfrac{\abs{\nabla z}^2}{z}
\\[3mm]\phantom{\Gd z----}
+(p-2)\myfrac{\langle \nabla z,\nabla u\rangle^2}{z^2}+
(q+2-p)z^{\frac{q-p}{2}}\langle\nabla z,\nabla u\rangle-(N-1)B^2z.
\EA\ee
Since $z^{\frac{q-p}{2}}\abs{\langle\nabla z,\nabla u\rangle}\leq z^{\frac{q+1-p}{2}}\myfrac{\abs{\nabla z}}{\sqrt z}$, we can take $a=a(p,q)>0$ large enough so that the right-hand side of (\ref{W3}) is bounded from below  by
$Cz^{q+2-p}-D\myfrac{\abs{\nabla z}^2}{z}$
for some $C,\,D>0$ which depend only on $p$ and $q$.
We set $$\CA (v):=-\Gd_2 v-(p-2)\myfrac{\langle D^2v(\nabla u),\nabla u\rangle}{\abs{\nabla u}^2}=-\sum_{i,j=1}^Na_{ij}v_{x_ix_j}$$
where the $a_{ij}$ depend on $\nabla u$ and satisfy
\begin{equation*}
\gth\abs\xi^2\leq \sum_{i,j=1}^Na_{ij}\xi_i\xi_j\leq \Gth\abs\xi^2\qquad\forall \xi=(\xi_1,...,\xi_n)\in\BBR^n,
\end{equation*}
where $\gth=\min\{1,p-1\}$ and  $\Gth=\max\{1,p-1\}$. Then
\bel{W5}
\CL^*(z):=\CA (z)+Cz^{q+2-p}- D\myfrac{\abs{\nabla z}^2}{z}-(n-1)B^2z\leq 0\qquad\text {in }M_+.
\ee

The next lemma is a local estimate. 
%%%%%%%%%%%%%%%%%%%%%%%%%%%%%%%%%%%%%%%%%%%%%%%%%%%%%%%%%%%%%%%%%%%%%%%%%%%%%%%%%%%%%%%%%%%%%%%%%%%%%%%%%%%%%%%%%%%%%%%%%%%LEMMA%%%%%%%%%%%%%%%%%%%%%%%%%%%%%%%%%%%%%%%%%%%%%%%%%%%%%%%%%%%%%%%%%%%%%%%%%%%%%%%%%%%%%%%%%%%%%
\blemma{compa} Let $B_R(a)\subset M^n(g)$ be the geodesic ball of radius $R>0$ and center $a$. Assume that
$Ricc_g\geq -(n-1)B^2$ and  either $Sec_g\geq -S^2$ for some $S^2:=S_R^2$ in $B_R(a)$ if $p>2$, or $Sec_g\leq 0$ if $1<p<2$. Then there exists $c=c(n,p,q)>0$ such that the function 
\bel{Ge0}w(x)=\gl\left(R^2-r^2(x)\right)^{-\frac{2}{q+1-p}}+\gm\quad \text {with }\,r=r(x)=d(x,a),
\ee
satisfies $\CL^*(w)\geq 0$ in $B_R(a)$, provided
\bel{Ge2}\BA{l}
\gl=c\max\left\{(R^4B^2)^{\frac{1}{q+1-p}},((1+B+(p-2)_+S)R^3)^{\frac{1}{q+1-p}}\right\}
\EA\ee
and
\bel{Ge2'}
\gm=((n-1)B^2)^{\frac{1}{q+1-p}}.
\ee
\es
\noindent {\bf Proof. } We recall that $ \Gd_2w=w''+w'\Gd_2r $ and by \cite[Lemma 1]{RRV} 
$$\Gd_2 r\leq (n-1)B\coth ( Br)\leq \myfrac{n-1}{r}\left(1+Br\right).$$ 
Then
\bel{Ge4}
\Gd_2w\leq \frac{4}{q+1-p}(R^2-r^2)^{-\frac{2(q+2-p)}{q+1-p}}\left(\frac{2r^2(q+3-p)}{q+1-p}+(R^2-r^2)(1+(n-1)\left(1+Br\right)\right).
\ee
Moverover from \cite[Chap 2, p. 23]{GW}
\bel{Ge5}
D^2w=w''dr\otimes dr+w'D^2r.
\ee
If $0\geq Sec_g(x)\geq -S^2$, there holds
\bel{Ge6}
0\leq D^2r\leq S\coth (Sr)g\leq \myfrac{S}{r}\left(1+Sr\right)g.
\ee
Therefore, if $p\geq 2$ and $Sec_g\geq -S^2$, we get

\bel{Ge7}
\myfrac{\langle D^2w(\nabla u),\nabla u\rangle}{\abs{\nabla u}^2}\leq 
\frac{4}{q+1-p}(R^2-r^2)^{-\frac{2(q+2-p)}{q+1-p}}\left(\frac{2r^2(q+3-p)}{q+1-p}+(R^2-r^2)(2+Sr)\right),
\ee
while, if $p\leq 2$ and $Sec_g\leq 0$,

\bel{Ge7'}\myfrac{\langle D^2w(\nabla u),\nabla u\rangle}{\abs{\nabla u}^2}\leq\frac{4}{q+1-p}(R^2-r^2)^{-\frac{2(q+2-p)}{q+1-p}}\left(\frac{2r^2(q+3-p)}{q+1-p}+2(R^2-r^2)\right).
\ee
As a consequence
\bel{Ge8}
\BA {lll}\CA(w)= -\Gd w-(p-2)\myfrac{\langle D^2w(\nabla u),\nabla u\rangle}{\abs{\nabla u}^2}\\[2mm]
\phantom{\CA(w)}
\geq  -k\gl(R^2-r^2)^{-\frac{2(q+2-p)}{q+1-p}}(R^2+(R^2-r^2)B_pr)
\EA\ee
for some $k=k(n,p,q)$, where $B_p=B+(p-2)_+S$. Since
$$w^{q+2-p}\geq \gl^{q+2-p}\left(R^2-r^2\right)^{-\frac{2(q+1-p)}{q+1-p}}+\gm^{q+2-p},
$$
we have
\bel{Ge9}\BA{l}
\CL^*(w)\geq \gl(R^2-r^2)^{-\frac{2(q+2-p)}{q+1-p}}\left(-k(R^2+(R^2-r^2)B_pr)-D\myfrac{16}{(q+1-p)^2}r^2
+C\gl^{q+1-p}\right)\\[4mm]\phantom{\CL^*(w)}
+\gm^{q+2-p}-(n-1)B^2\gl\left(R^2-r^2\right)^{-\frac{2}{q+1-p}}-(n-1)B^2\gm.
\EA\ee
We first take 
\bel{Ge10}\gm=((n-1)B^2)^{\frac{1}{q+1-p}}.
\ee
Next we choose $\gl$ in order to have, uniformly for $0\leq r<R$,
$$
2^{-1}C\gl^{q+1-p}\geq k\left(R^2+(R^2-r^2)B_pr\right)+\myfrac{16Dr^2}{(q+1-p)^2}
$$
and
$$2^{-1}C\gl^{q+2-p}(R^2-r^2)^{-\frac{2(q+2-p)}{q+1-p}}\geq (n-1)B^2\gl\left(R^2-r^2\right)^{-\frac{2}{q+1-p}}.
$$
There exists $c=c(n,p,q)$ such that, if
\bel{Ge11}
\gl=c\max\left\{(R^4B^2)^{\frac{1}{q+1-p}},((1+B_p)R^3)^{\frac{1}{q+1-p}}\right\},
\ee
then $\CL^*(w)\geq 0$ holds.
%%%%%%%%%%%%%%%%%%%%%%%%%%%%%%%%%%%%%%%%%%%%%%%%%%%%%%%%%%%%%%%%%%%%%%%%%%%%%%%%%%%%%%%%%%%%%%%%%%%%%%%%%%%%%%%%%%%%%%%%%%%LEMMA%2-%%%%%%%%%%%%%%%%%%%%%%%%%%%%%%%%%%%%%%%%%%%%%%%%%%%%%%%%%%%%%%%%%%%%%%%%%%%%%%%%%%%%%%%%%%%%%

\blemma{compa2} Under the assumptions of \rlemma{compa}, any $C^1$ solution of (\ref{equ1-2}) in $M$ satisfies
\bel{est1}
\abs{\nabla u(x)}\leq c_{n,p,q}\max\left\{B^{\frac{1}{q+1-p}},(1+B_p)^{\frac{1}{2(q+1-p)}}(d(x,\prt\Gw))^{-\frac{1}{2(q+1-p)}}\right\}\qquad\forall x\in \Gw,
\ee 
for every  domain $\Gw\subset M$, where $B_p=B+(p-2)_+S$ and $S=S_{d(x,\prt\Gw)}$.
\es
\noindent {\bf Proof. } Assume $a\in \Gw$, with $R<d(a,\prt\Gw)$. Let $w$ be as in \rlemma{compa}, then in any  connected component  $G$ of $\{x\in B_R(a):z(x)-w(x)>0\}$ we find
\bel{Ge12}
\CA (z-w)+C\left(z^{q+2-p}-w^{q+2-p}\right)-(n-1)B^2(z-w) -D\left(\myfrac{\abs{\nabla z}^2}{z}-\myfrac{\abs{\nabla w}^2}{w}\right)\leq 0.
\ee
By the mean value theorem and since $w(a)$ is the minimum of $w$, there holds 
\bel{Ge13}
C(z^{q+2-p}-w^{q+2-p})-(n-1)B^2(z-w)>0,
\ee
provided $C(q+2-p)(w(a))^{q+1-p}>(n-1)B^2$. Since  $w(a)>\mu=\left((n-1)B^2\right)^{\frac{1}{q+1-p}}$ and $q+2-p>1$, this condition is fulfilled, up to replacing $\mu$ by $A\mu$ for some $A=A(p,q)>1$. If $x_0\in G$ is such that 
$z-w$ is maximal at $x_0$, we derive that 
$$\CA (z-w)+C\left(z^{q+2-p}-w^{q+2-p}\right)-(N-1)B^2(z-w) -D\left(\myfrac{\abs{\nabla z}^2}{z}-\myfrac{\abs{\nabla w}^2}{w}\right)\leq 0
$$
if $x=x_0$, which is a contradiction. Thus $G=\emptyset$, $z\leq w$  and (\ref{est1}) follows.\medskip

The proof of Theorem 1 and Corollary 2 follows by taking $\Omega=B_R(x)$ and letting $R\to\infty$.\medskip

\noindent{\bf Proof of Theorem 3.} We take $q=p$ and assume that $v$ is $p$-harmonic and positive. If we write $v=e^{-\frac{u}{p-1}}$ , then $u$ satisfies 
$$-\Gd_p u+\abs{\nabla u}^p=0.$$
If $Ricc_g(x)\geq 0$, $u$ is constant by Corollary 2, and so is $v$. If $\inf\{Ricc_g(x):x\in M\}=(1-n)B^2<0$
we apply $(\ref{est1})$ to $\nabla u$. If $\gamma$ is a minimizing geodesic from $a$ to $x$, then $\abs{\gamma'(t)}=1$
and 
$$u(x)-u(a)=\int_{0}^{d(x,a)}\myfrac{d}{dt}u\circ\gamma (t) dt=\int_{0}^{d(x,a)}\langle\nabla u\circ\gamma (t),\gamma' (t)\rangle dt.
$$
Since
$$\abs{\langle\nabla u\circ\gamma (t),\gamma' (t)}\leq \abs{ \nabla u\circ\gamma (t)}\leq c_{n,p,p}B,
$$
we obtain
\bel{est3}
u(a)-c_{n,p,p}B\dist(x,a)\leq u(x)\leq u(a)+c_{n,p,p}B\dist(x,a)\qquad\forall x\in M.
\ee
Then (\ref{equ4'-2}) follows since $u=(1-p)\ln v$.

%%%%%%%%%%%%%%%%%%%

%%%%%bibliographie%%%%%%%%%%%%%%%%%%%%%%%%%%%%%%%%%%
%%%%%%%%%%%%%%%%%%% %%%%%%%%%%%%%%%%%%%%%%%%%%%%

\end {document}